\documentclass[final,natbib]{svjour3}                     
\smartqed  
\usepackage{graphicx,amsmath,amssymb,siunitx,soul}
\usepackage[table]{xcolor}
\usepackage[colorlinks=true,breaklinks=true,citecolor=blue,linkcolor=blue,urlcolor=blue,allcolors=blue]{hyperref}
\usepackage{arydshln}
\definecolor{tableshade}{RGB}{198.6000  227.1000  239.1000}
\def\figw{0.49\textwidth}

\begin{document}

\title{New semi-analytical solutions for advection-dispersion equations in multilayer porous media
}

\titlerunning{New semi-analytical solutions for advection-dispersion equations in multilayer porous media}

\author{Elliot J Carr}


\institute{Elliot J Carr \at
              School of Mathematical Sciences, Queensland University of Technology (QUT), Brisbane, Australia. \\
              \email{elliot.carr@qut.edu.au}
}

\date{}

\maketitle

\begin{abstract}
A new semi-analytical solution to the advection-dispersion-reaction equation for modelling solute transport in layered porous media is derived using the Laplace transform. Our solution approach involves introducing unknown functions representing the dispersive flux at the interfaces between adjacent layers, allowing the multilayer problem to be solved separately on each layer in the Laplace domain before being numerically inverted back to the time domain. The derived solution is applicable to the most general form of linear advection-dispersion-reaction equation, a finite medium comprising an arbitrary number of layers, continuity of concentration and dispersive flux at the interfaces between adjacent layers and transient boundary conditions of arbitrary type at the inlet and outlet. The derived semi-analytical solution extends and addresses deficiencies of existing analytical solutions in a layered medium, which consider analogous processes such as diffusion or reaction-diffusion only and/or require the solution of complicated nonlinear transcendental equations to evaluate the solution expressions. Code implementing our semi-analytical solution is supplied and applied to a selection of test cases, with the reported results in excellent agreement with a standard numerical solution and other analytical results available in the literature.
\end{abstract}

\section{Introduction}
\begin{figure*}[t]
\centering
\fbox{\includegraphics[width=0.9\textwidth]{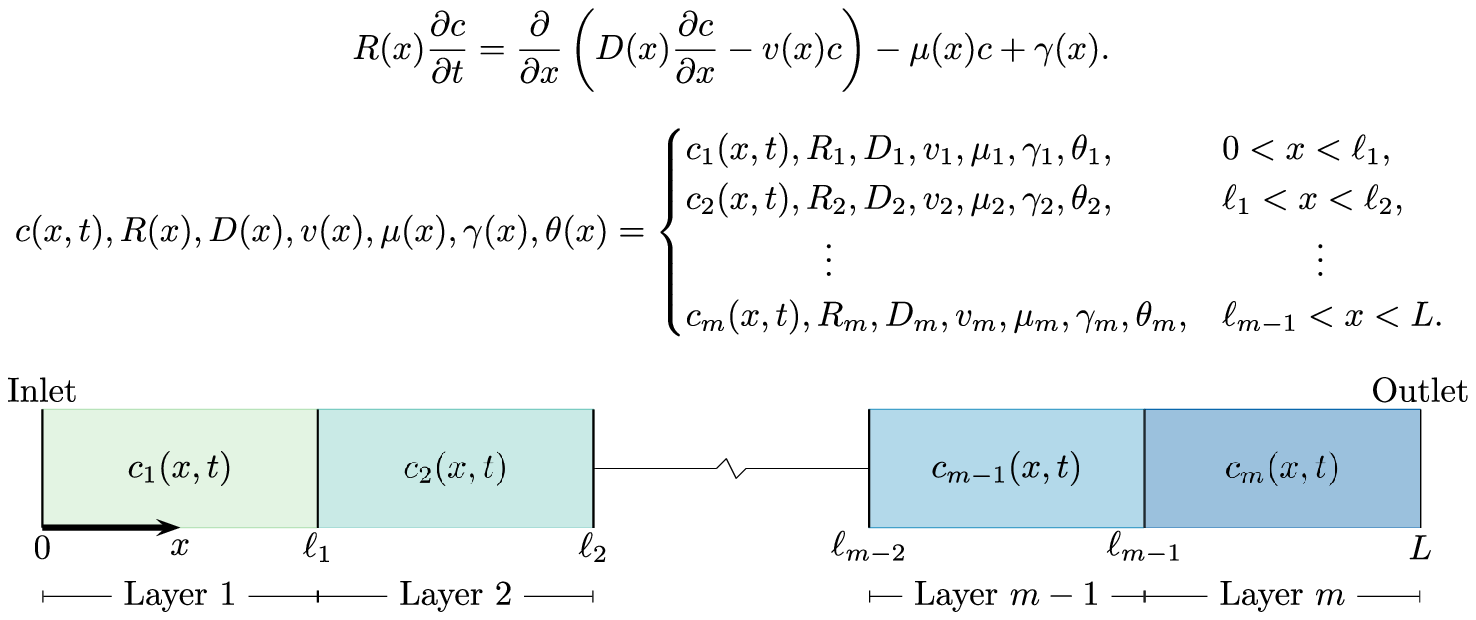}}
\caption{Advection-dispersion-reaction in an $m$-layered medium. We solve for the concentration $c(x,t)$, where the retardation factor $R$, dispersion coefficient $D$, pore-water velocity $v$, volumetric water content $\theta$ and rate constants for first-order decay $\mu$ and zero-order production $\gamma$ are constant within each layer but vary across layers. Continuity of concentration and dispersive flux are imposed at the interfaces between adjacent layers ($x = \ell_{i}$, $i = 1,\hdots,m-1$) and transient boundary conditions of arbitrary-type are specified at the inlet ($x=0$) and outlet ($x=L$).}
\label{fig:Figure1}
\end{figure*}

Solute transport in porous media due to dispersion and groundwater flow is typically modelled using the advection-dispersion-reaction equation \citep{van_genuchten_1982,leij_1991,liu_1998,goltz_2017}. While solving this equation in heterogeneous porous media usually requires the application of numerical methods, analytical solutions are generally preferred when available as they are exact and continuous in space and time. The focus of this paper is analytical solutions for layered porous media, which are commonly observed in natural and constructed environments such as stratified soils and landfill clay liners \citep{liu_1998}. For this problem the transport coefficients (dispersion coefficient, retardation factor etc) are piecewise constant with the goal being to obtain the solute concentration in each layer as a function of space and time (see Figure \ref{fig:Figure1}).

Analytical solutions of advection-dispersion-reaction equations (and closely related equations such as advection-diffusion and reaction-diffusion equations) in layered media continue to attract interest \citep{zimmerman_2016,yang_2017,guerrero_2013,carr_2018b,rodrigo_2016}. Mainly due to their prevalence in heat conduction problems, most analytical solutions in the literature are developed for multilayer diffusion problems \citep{carr_2016,carr_2018b,sun_2004,monte_2002,hickson_2009,rodrigo_2016} with significantly less literature concerning advection-dispersion, reaction-diffusion or advection-dispersion-reaction problems. 

Analytical solutions for advection-dispersion equations in layered media have been presented by a limited number of authors. Amongst the earliest work in this area is that of \citet{leij_1991}, who applied the Laplace transform to solve the advection-dispersion equation (with retardation factor) on a semi-infinite two-layer medium with finite first layer, semi-infinite second layer and continuity of concentration and dispersive flux at the interfaces between adjacent layers. Both concentration-type and flux-type boundary conditions were considered at the inlet and a zero concentration gradient was applied at the outlet (Figure \ref{fig:Figure1}). Exact expressions for the concentration in the Laplace domain were obtained and numerically inverted.

In follow up work, \citet{leij_1995} derived approximate analytical solutions by first expanding the Laplace domain concentration in both layers as infinite series, truncating each series after the first term and employing analytical inversion of the Laplace transform to convert the concentration back to the time domain. Subsequent analytical solutions for finite layered media and an arbitrary number of layers were derived by \citet{liu_1998} and then later by \citet{guerrero_2013}, both using the method of eigenfunction expansion. \citet{liu_1998} considered the advection-dispersion equation and \citet{guerrero_2013} the advection-diffusion-reaction equation with a first order decay term. Both papers treated inlet boundary conditions of flux-type only with \citet{liu_1998} allowing for an arbitrary time-varying inlet concentration and \citet{guerrero_2013} a constant inlet concentration. Both approaches require a transcendental equation to be solved numerically for the eigenvalues with the accuracy of the solutions depending on the number of eigenvalues used in the expansions. 

Recent work has revisited the use of Laplace transforms for solving multilayer transport problems due to the ease of treating different boundary conditions and the ability to avoid solving complicated nonlinear transcendental equations as required by eigenfunction expansion solutions \citep{carr_2016}.  Namely, \citet{carr_2016} and \citet{rodrigo_2016} solved the multilayer diffusion problem using the Laplace transform for an arbitrary number of layers and various types of boundary and interface conditions. \citet{zimmerman_2016} solved a reaction-diffusion equation (without retardation factor) with a first order reaction term on a finite medium consisting of an arbitrary number of layers. All three of these papers \citep{carr_2016,rodrigo_2016,zimmerman_2016} introduce an unknown function of time at each interface, representing either the diffusive/dispersive flux or concentration, which allows the multilayer problem to be isolated and solved separately on each layer in the Laplace domain before being inverted back to the time domain. 

Further relevant work involving the Laplace transform has focussed on the permeable reaction barrier-aquifer problem. In particular, semi-analytical solutions for the advection-dispersion-reaction equation have been developed by {\cite{park_2009}} for a semi-infinite one-dimensional medium using the Laplace transform and {\cite{chen_2016}} for a finite two-dimensional medium using Laplace and finite Fourier transforms. Both solutions are limited to two layers and numerically invert the Laplace transform using the {\cite{dehoog_1982}} algorithm.

The aim of the current paper is to extend, generalize and merge the work of \citet{carr_2016}, \citet{rodrigo_2016}, \citet{zimmerman_2016} and \citet{park_2009} to solve advection-dispersion-reaction problems in one-dimensional media comprising an arbitrary number of layers using the Laplace transform. The derived solution is applicable to the most general form of linear advection-dispersion-reaction equation \citep{van_genuchten_1982}, a finite medium comprising an arbitrary number of layers, continuity of concentration and dispersive flux at the interfaces between adjacent layers and transient boundary conditions of arbitrary type at the inlet and outlet. To the best of our knowledge, analytical solutions satisfying all the above conditions have not previously appeared in the published literature. The derived solutions and resulting code supersedes our previous work on multilayer diffusion \citep{carr_2016,carr_2018a} by removing the requirement of numerically computing eigenvalues in each layer. 

The remaining sections of this paper are organised as follows. In section \ref{sec:transport_model} we describe the multilayer advection-dispersion-reaction problem considered in this work. Section \ref{sec:solution} develops our proposed solution procedure using the Laplace transform. In section \ref{sec:results}, the developed solutions are applied to a wide selection of test cases and compared to other analytical and numerical solutions available in the literature. In section \ref{sec:conclusions}, we summarise the work and discuss possible avenues for future research.

\section{Multilayer transport model}
\label{sec:transport_model}
We consider solute transport across an $m$-layered porous medium partitioned as $0 = \ell_{0} < \ell_{1} < \cdots < \ell_{m-1} < \ell_{m} = L$ (Figure \ref{fig:Figure1}). Let $c_{i}(x,t)$ be the solute concentration [$\text{M}\text{L}^{-3}$] in the $i$th layer, where $x\in (\ell_{i-1},\ell_{i})$ is the distance from the inlet at $x = 0$ and $t > 0$ is time. The governing transport equation in the $i$th layer is 
\begin{gather}
\label{eq:pde}
R_{i}\frac{\partial c_{i}}{\partial t} = D_{i}\frac{\partial^{2}c_{i}}{\partial x^{2}} - v_{i}\frac{\partial c_{i}}{\partial x} - \mu_{i}c_{i} + \gamma_{i},
\end{gather}
where $R_{i}>0$ is the retardation factor $[-]$, $D_{i}>0$ is the dispersion coefficient [$\text{L}^{2}\text{T}^{-1}$], $v_{i}$ is the pore-water velocity [$\text{L}\text{T}^{-1}$], $\mu_{i}$ is the rate constant for first-order decay [$\text{M}\text{L}^{-3}\text{T}^{-1}$] and $\gamma_{i}$ is the rate constant for zero-order production [$\text{T}^{-1}$] \citep{van_genuchten_1982}. Throughout this paper we will denote by $c(x,t)$ the result of amalgamating the layer concentrations, $c_{i}(x,t)$ ($i = 1,\hdots,m$), as defined in Figure \ref{fig:Figure1}.

The transport equations (\ref{eq:pde}) are accompanied by initial, interface and boundary conditions. The concentration is initially assumed constant in each layer
\begin{gather}
\label{eq:ic}
c_{i}(x,0) = f_{i},
\end{gather}
concentration and dispersive flux is assumed continuous at the interfaces between adjacent layers ($i = 1,\hdots,m-1$) \citep{leij_1991,liu_1998,guerrero_2013}:
\begin{gather}
\label{eq:concentration_continuity}
c_{i}(\ell_{i},t) = c_{i+1}(\ell_{i},t),\\
\label{eq:dispersive_flux_continuity}
\theta_{i}D_{i}\frac{\partial c_{i}}{\partial x}(\ell_{i},t) = \theta_{i+1}D_{i+1}\frac{\partial c_{i+1}}{\partial x}(\ell_{i},t),
\end{gather}
and general Robin boundary conditions are considered at the inlet ($x = 0$) and outlet ($x = L$):
\begin{gather}
\label{eq:general_bcL}
a_{0}c_{1}(0,t) - b_{0}\frac{\partial c_{1}}{\partial x}(0,t) = g_{0}(t),\\
\label{eq:general_bcR}
a_{L}c_{m}(L,t) + b_{L}\frac{\partial c_{m}}{\partial x}(L,t) = g_{L}(t).
\end{gather}
The constant $\theta_{i}$ appearing in the interface condition (\ref{eq:dispersive_flux_continuity}) is the volumetric water content [$\text{L}^{3}\text{L}^{-3}$] in the $i$th layer \citep{van_genuchten_1982}. In the boundary conditions (\ref{eq:general_bcL})--(\ref{eq:general_bcR}), $a_{0}$, $b_{0}$, $a_{L}$ and $b_{L}$ are constants and $g_{0}(t)$ and $g_{L}(t)$ are arbitrary specified functions of time with the subscripts $0$ and $L$ denoting the inlet ($x = 0$) and outlet ($x = L$), respectively (Figure \ref{fig:Figure1}). We remark that $b_{0}$ and $b_{L}$ are both non-negative, at least one of $a_{0}$ or $b_{0}$ must be non-zero and at least one of $a_{L}$ and $b_{L}$ must be non-zero. 

For solute transport problems, commonly \citep{leij_1991,liu_1998,goltz_2017} either a concentration-type boundary condition
\begin{gather}
\label{eq:inlet_concentration_type}
c_{1}(0,t) = c_{0}(t),
\end{gather}
or a flux-type boundary condition
\begin{gather}
\label{eq:inlet_flux_type}
v_{1}c_{1}(0,t) - D_{1}\frac{\partial c_{1}}{\partial x}(0,t) = v_{1}c_{0}(t),
\end{gather}
is applied at the inlet, where $c_{0}(t)$ is the specified inlet concentration, while a zero concentration gradient is applied at the outlet:
\begin{gather}
\label{eq:outlet_zero_concentration_gradient}
\frac{\partial c_{m}}{\partial x}(L,t) = 0.
\end{gather}
These boundary conditions are obtained from the general boundary conditions (\ref{eq:general_bcL})--(\ref{eq:general_bcR}) by setting $a_{0} = 1$, $b_{0} = 0$, $g_{0}(t) = c_{0}(t)$ [Eq (\ref{eq:inlet_concentration_type})], $a_{0} = v_{1}$, $b_{0} = D_{1}$, $g_{0}(t) = v_{1}c_{0}(t)$ [Eq (\ref{eq:inlet_flux_type})] and $a_{L} = 0$, $b_{L} = 1$, $g_{L}(t) = 0$ [Eq (\ref{eq:outlet_zero_concentration_gradient})], respectively. 

\section{Semi-analytical solution}
\label{sec:solution}

\subsection{General solution for $m$ layers in Laplace space}
\label{sec:solution1}
To solve the multilayer transport model (\ref{eq:pde})--(\ref{eq:general_bcR}), we reformulate the model into $m$ isolated single layer problems \citep{carr_2016,rodrigo_2016,zimmerman_2016}. Introducing unknown functions of time, $g_{i}(t)$ ($i = 1,\hdots,m-1$), to denote the following scalar multiple of the (negative) dispersive flux at the layer interfaces \citep{carr_2016,rodrigo_2016}:
\begin{align*}
g_{i}(t) = \theta_{i}D_{i}\frac{\partial c_{i}}{\partial x}(\ell_{i},t),
\end{align*}
yields the following equivalent form for the multilayer transport model (\ref{eq:pde})--(\ref{eq:general_bcR}):

\bigskip\noindent \textbf{First layer} ($i = 1$)
\begin{gather}
\label{eq:pde1}
R_{1}\frac{\partial c_{1}}{\partial t} = D_{1}\frac{\partial^{2}c_{1}}{\partial x^{2}} - v_{1}\frac{\partial c_{1}}{\partial x} - \mu_{1}c_{1} + \gamma_{1},\\
\label{eq:ic1}
c_{1}(x,0) = f_{1},\\
\label{eq:bcL1}
a_{0}c_{1}(0,t) - b_{0}\frac{\partial c_{1}}{\partial t}(0,t) = g_{0}(t),\\
\label{eq:bcR1}
\theta_{1}D_{1}\frac{\partial c_{1}}{\partial x}(\ell_{1},t) = g_{1}(t),
\end{gather}

\medskip\noindent \textbf{Middle layers} ($i = 2,\hdots,m-1$)
\begin{gather}
\label{eq:pdei}
R_{i}\frac{\partial c_{i}}{\partial t} = D_{i}\frac{\partial^{2}c_{i}}{\partial x^{2}} - v_{i}\frac{\partial c_{i}}{\partial x} - \mu_{i}c_{i} + \gamma_{i},\\
\label{eq:ici}
c_{i}(x,0) = f_{i},\\
\label{eq:bcLi}
\theta_{i}D_{i}\frac{\partial c_{i}}{\partial x}(\ell_{i-1},t) = g_{i-1}(t),\\
\label{eq:bcRi}
\theta_{i}D_{i}\frac{\partial c_{i}}{\partial x}(\ell_{i},t) = g_{i}(t),
\end{gather}

\medskip\noindent \textbf{Last layer} ($i = m$)
\begin{gather}
\label{eq:pdem}
R_{m}\frac{\partial c_{m}}{\partial t} = D_{m}\frac{\partial^{2}c_{m}}{\partial x^{2}} - v_{m}\frac{\partial c_{m}}{\partial x} - \mu_{m}c_{m} + \gamma_{m},\\
\label{eq:icm}
c_{m}(x,0) = f_{m},\\
\label{eq:bcLm}
\theta_{m}D_{m}\frac{\partial c_{m}}{\partial x}(\ell_{m-1},t) = g_{m-1}(t),\\
\label{eq:bcRm}
a_{L}c_{m}(L,t) + b_{L}\frac{\partial c_{m}}{\partial x}(L,t) = g_{L}(t),
\end{gather}
with each problem coupled together by imposing continuity of concentration at the interfaces between adjacent layers (\ref{eq:concentration_continuity}) \citep{carr_2016,rodrigo_2016,carr_2018a}. 

We remark that the advection-dispersion-reaction equations (\ref{eq:pde1}), (\ref{eq:pdei}) and (\ref{eq:pdem}) can be reduced to standard heat/diffusion equations via the change of variables: 
\begin{align}
\label{eq:sh_transformation}
c_{i}(x,t) = \frac{\gamma_{i}}{\mu_{i}} + \exp\left(\frac{v_{i}x}{2D_{i}} - \left[\frac{v_{i}^{2}}{4R_{i}D_{i}} + \frac{\mu_{i}}{R_{i}}\right]t\right)u_{i}(x,t). 
\end{align} 
This seems to suggest that the transformed problem is amendable to solution by methods designed for the multilayer diffusion problem presented in previous work (e.g. \cite{carr_2016}, etc). However, this is not the case as the transformation (\ref{eq:sh_transformation}) gives rise to an advection term in the interface conditions (\ref{eq:bcR1}), (\ref{eq:bcLi})--(\ref{eq:bcRi}) and (\ref{eq:bcLm}) for the transformed variable $u_{i}(x,t)$ that would require the solutions to be redeveloped. Moreover, the transformation complicates the problem converting the initial conditions (\ref{eq:ic1}), (\ref{eq:ici}), (\ref{eq:icm}) from constants to spatially-dependent. We therefore do not invoke the transformation (\ref{eq:sh_transformation}) and take an alternative approach in this paper.

To solve the isolated single-layer problems (\ref{eq:pde1})--(\ref{eq:bcR1}), (\ref{eq:pdei})--(\ref{eq:bcRi}) and (\ref{eq:pdem})--(\ref{eq:bcRm}), we take Laplace transforms yielding the boundary value problems:

\bigskip\noindent \textbf{First layer} ($i = 1$)
\begin{gather}
\label{eq:bvp1_de}
D_{1}C_{1}'' - v_{1}C_{1}' - (\mu_{1} + R_{1}s)C_{1} = - R_{1}f_{1} - \frac{\gamma_{1}}{s},\\
a_{0}C_{1}(0,s) - b_{0}C_{1}'(0,s) = G_{0}(s),\\
\label{eq:bvp1_bc2}
\theta_{1}D_{1}C_{1}'(\ell_{1},s) = G_{1}(s),
\end{gather}

\medskip\noindent \textbf{Middle layers} ($i = 2,\hdots,m-1$)
\begin{gather}
\label{eq:bvpi_de}
D_{i}C_{i}'' - v_{i}C_{i}' - (\mu_{i} + R_{i}s)C_{i} = - R_{i}f_{i} - \frac{\gamma_{i}}{s},\\
\theta_{i}D_{i}C_{i}'(\ell_{i-1},s) = G_{i-1}(s),\\
\label{eq:bvpi_bc2}
\theta_{i}D_{i}C_{i}'(\ell_{i},s) = G_{i}(s),
\end{gather}

\medskip\noindent \textbf{Last layer} ($i = m$)
\begin{gather}
\label{eq:bvpm_de}
D_{m}C_{m}'' - v_{m}C_{m}' - (\mu_{m} + R_{m}s)C_{m} = - R_{m}f_{m} - \frac{\gamma_{m}}{s},\hspace{-1cm}\\
\theta_{m}D_{m}C_{m}'(\ell_{m-1},s) = G_{m-1}(s),\\
\label{eq:bvpm_bc2}
a_{L}C_{m}(L,s) + b_{L}C_{m}'(L,s) = G_{L}(s),
\end{gather}
where the prime notation ($'$) denotes a derivative with respect to $x$, $C_{i}(x,s) = \mathcal{L}\{c_{i}(x,t)\}$ denotes the Laplace transform of $c_{i}(x,t)$ with transformation variable $s\in\mathbb{C}$ and $G_{i}(s) = \mathcal{L}\{g_{i}(t)\}$ for $i = 1,\hdots,m-1$. Both $G_{0}(s) = \mathcal{L}\{g_{0}(t)\}$ and $G_{L}(s) = \mathcal{L}\{g_{L}(t)\}$ are assumed to be able to be found analytically. The boundary value problems (\ref{eq:bvp1_de})--(\ref{eq:bvp1_bc2}), (\ref{eq:bvpi_de})--(\ref{eq:bvpi_bc2}) and (\ref{eq:bvpm_de})--(\ref{eq:bvpm_bc2}) all involve second-order constant-coefficient differential equations, which can be solved using standard techniques to give the following expressions for the concentration in the Laplace domain:
\begin{alignat}{2}
\label{eq:C1}
C_{1}(x,s) &= A_{1}(x,s)G_{0}(s) + B_{1}(x,s)G_{1}(s) + P_{1}(x,s),\\
\label{eq:Ci}
C_{i}(x,s) &= A_{i}(x,s)G_{i-1}(s) + B_{i}(x,s)G_{i}(s) + P_{i}(x,s),\quad i = 2,\hdots,m-1,\\
\label{eq:Cm}
C_{m}(x,s) &= A_{m}(x,s)G_{m-1}(s) + B_{m}(x,s)G_{L}(s) + P_{m}(x,s),
\end{alignat}
where the functions $P_{i}$, $A_{i}$ and $B_{i}$ ($i = 1,\hdots,m$) are defined in Table \ref{tab:FAB}.

\begin{table*}
\centering
\caption{Definition of the functions $P_{i}(x,s)$, $A_{i}(x,s)$ and $B_{i}(x,s)$ ($i = 1,\hdots,m$) appearing in the Laplace domain concentrations (\ref{eq:C1})--(\ref{eq:Cm}).}
\renewcommand{\arraystretch}{1.6}
\begin{tabular}{l}
\hline
\textbf{First layer} ($i = 1$)\\ 
$\displaystyle P_{1}(x,s) = \Psi_{1}(s) + \frac{a_{0}}{\beta_{1}(s)}\Bigl\{\lambda_{1,1}(s)\Psi_{1,2}(x,s) - \lambda_{1,2}(s)\Psi_{1,2}(\ell_{1},s)\Psi_{1,1}(x,s)\Bigr\}\Psi_{1}(s)$,\\
$\displaystyle A_{1}(x,s) = \frac{1}{\beta_{1}(s)}\Bigl\{\lambda_{1,2}(s)\Psi_{1,2}(\ell_{1},s)\Psi_{1,1}(x,s) - \lambda_{1,1}(s)\Psi_{1,2}(x,s)\Bigr\},$\\
$\displaystyle B_{1}(x,s) = \frac{1}{\theta_{1}D_{1}\beta_{1}(s)}\Bigl\{\bigl[a_{0}-b_{0}\lambda_{1,1}(s)\bigr]\Psi_{1,1}(0,s)\Psi_{1,2}(x,s) - \bigl[a_{0} - b_{0}\lambda_{1,2}(s)\bigr]\Psi_{1,1}(x,s)\Bigr\}$,\\
\text{where}\\
$\displaystyle \beta_{1}(s) = [a_{0} - b_{0}\lambda_{1,1}(s)]\lambda_{1,2}(s)\exp(-[\lambda_{1,1}(s)-\lambda_{1,2}(s)]\ell_{1}) - [a_{0} - b_{0}\lambda_{1,2}(s)]\lambda_{1,1}(s)$.\\
\hline
\textbf{Middle layers} ($i = 2,\hdots,m-1$)\\
$\displaystyle P_{i}(x,s) = \Psi_{i}(s)$,\\
$\displaystyle A_{i}(x,s) = \frac{1}{\theta_{i}D_{i}\beta_{i}(s)}\Bigl\{\lambda_{i,2}(s)\Psi_{i,2}(\ell_{i},s)\Psi_{i,1}(x,s) - \lambda_{i,1}(s)\Psi_{i,2}(x,s)\Bigr]$,\\
$\displaystyle B_{i}(x,s) = \frac{1}{\theta_{i}D_{i}\beta_{i}(s)}\Bigl\{\lambda_{i,1}(s)\Psi_{i,1}(\ell_{i-1},s)\Psi_{i,2}(x,s)-\lambda_{i,2}(s)\Psi_{i,1}(x,s)\Bigr\}$,\\
\text{where}\\
$\displaystyle \beta_{i}(s) = \lambda_{i,1}(s)\lambda_{i,2}(s)\left\{\exp(-[\lambda_{i,1}(s)-\lambda_{i,2}(s)](\ell_{i}-\ell_{i-1}))-1\right\}$.\\
\hline
\textbf{Last layer} ($i = m$)\\
$\displaystyle P_{m}(x,s) = \Psi_{m}(s) + \frac{a_{L}}{\beta_{m}(s)}\Bigl\{\lambda_{m,2}(s)\Psi_{m,1}(x,s) - \lambda_{m,1}(s)\Psi_{m,1}(\ell_{m-1},s)\Psi_{m,2}(x,s)\Bigr\}\Psi_{m}(s)$,\\
$\displaystyle A_{m}(x,s) = \frac{1}{\theta_{m}D_{m}\beta_{m}(s)}\Bigl\{\bigl[a_{L}+b_{L}\lambda_{m,2}(s)\bigr]\Psi_{m,2}(L,s)\Psi_{m,1}(x,s) - \left[a_{L}+b_{L}\lambda_{m,1}(s)\right]\Psi_{m,2}(x,s)\Bigr\}$,\\
$\displaystyle B_{m}(x,s) = \frac{1}{\beta_{m}(s)}\Bigl\{\lambda_{m,1}(s)\Psi_{m,1}(\ell_{m-1},s)\Psi_{m,2}(x,s) - \lambda_{m,2}(s)\Psi_{m,1}(x,s)\Bigr\}$,\\
\text{where}\\
$\displaystyle \beta_{m}(s) = [a_{L}+b_{L}\lambda_{m,2}(s)]\lambda_{m,1}(s)\exp(-[\lambda_{m,1}(s)-\lambda_{m,2}(s)](\ell_{m}-\ell_{m-1})) - [a_{L}+b_{L}\lambda_{m,1}(s)]\lambda_{m,2}(s)$.\\
\hline
For all layers ($i = 1,\hdots,m$):\\
$\displaystyle\Psi_{i}(s) = \frac{\gamma_{i}s^{-1} + R_{i}f_{i}}{\mu_{i} + R_{i}s}$,\\
$\displaystyle\Psi_{i,1}(x,s) = \exp\left[\lambda_{i,1}(s)(x-\ell_{i})\right]$ recalling $\ell_{m} = L$,\\
$\displaystyle\Psi_{i,2}(x,s) = \exp\left[\lambda_{i,2}(s)(x-\ell_{i-1})\right]$ recalling $\ell_{0} = 0$,\\
$\displaystyle\lambda_{i,1}(s) = \frac{v_{i} + \sqrt{v_{i}^{2} + 4D_{i}(R_{i}s + \mu_{i})}}{2D_{i}}$,\\
$\displaystyle\lambda_{i,2}(s) = \frac{v_{i} - \sqrt{v_{i}^{2} + 4D_{i}(R_{i}s + \mu_{i})}}{2D_{i}}$.\\
\hline
\end{tabular}
\label{tab:FAB}
\end{table*} 

With all other variables in the expressions (\ref{eq:C1})--(\ref{eq:Cm}) defined, to determine $G_{1}(s),\hdots,G_{m-1}(s)$, the Laplace transformations of the unknown interface functions $g_{1}(t),\hdots,g_{m-1}(t)$, we enforce continuity of concentration (\ref{eq:concentration_continuity}) at each interface in the Laplace domain \citep{carr_2016,rodrigo_2016,carr_2018a}:
\begin{align}
\label{eq:concentration_continuity_Laplace}
C_{i}(\ell_{i},s) = C_{i+1}(\ell_{i},s),\quad i = 1,\hdots,m-1.
\end{align}
Substituting (\ref{eq:C1})--(\ref{eq:Cm}) into the system of equations (\ref{eq:concentration_continuity_Laplace}) yields a linear system for $\mathbf{x} = \left[G_{1}(s),\hdots,G_{m-1}(s)\right]^{T}$, expressible in matrix form as
\begin{align}
\label{eq:linear_system}
\mathbf{A}\mathbf{x} = \mathbf{b},
\end{align}
where $\mathbf{A} = (a_{i,j})\in\mathbb{C}^{(m-1)\times(m-1)}$ is a tridiagonal matrix and $\mathbf{b} = (b_{i})\in\mathbb{C}^{(m-1)}$ is a vector with entries:
\begin{align*}
a_{1,1} &= B_{1}(\ell_{1},s) - A_{2}(\ell_{1},s),\\
a_{1,2} &= -B_{2}(\ell_{1},s),\\
a_{i,i-1} &= A_{i}(\ell_{i},s),\quad i = 2,\hdots,m-2,\\
a_{i,i} &= B_{i}(\ell_{i},s) - A_{i+1}(\ell_{i},s),\quad i = 2,\hdots,m-2,\\
a_{i,i+1} &= -B_{i+1}(\ell_{i},s),\quad i = 2,\hdots,m-2,\\
a_{m-1,m-2} &= A_{m-1}(\ell_{m-1},s),\\
a_{m-1,m-1} &= B_{m-1}(\ell_{m-1},s) - A_{m}(\ell_{m-1},s),\\
b_{1} &= P_{2}(\ell_{1},s) - P_{1}(\ell_{1},s) - A_{1}(\ell_{1},s)G_{0}(s),\\
b_{i} &= P_{i+1}(\ell_{i},s) - P_{i}(\ell_{i},s),\quad i = 2,\hdots,m-2,\\
b_{m-1} &= P_{m}(\ell_{m-1},s) - P_{m-1}(\ell_{m-1},s) + B_{m}(\ell_{m-1},s)G_{L}(s).
\end{align*}
Solving the linear system (\ref{eq:linear_system}) allows the functions $G_{1}(s),\hdots,G_{m-1}(s)$ to be computed and hence the Laplace transform of the concentration (\ref{eq:C1})--(\ref{eq:Cm}) can be evaluated at any $x$ and $s$ in the Laplace domain.

\subsection{Simplification for two layers in Laplace space}
\label{sec:solution2}
The formulation presented in the previous section breaks down for $m = 2$ layers due to the assumption of middle layers. While the expressions for $C_{1}(x,s)$ (\ref{eq:C1}) and $C_{2}(x,s)$ (\ref{eq:Cm}) remain valid:
\begin{align*}
C_{1}(x,s) &= A_{1}(x,s)G_{0}(s) + B_{1}(x,s)G_{1}(s)\nonumber + P_{1}(x,s),\\
C_{2}(x,s) &= A_{2}(x,s)G_{1}(s) + B_{2}(x,s)G_{L}(s)\nonumber + P_{2}(x,s),
\end{align*}
the defined entries of the linear system (\ref{eq:linear_system}) are no longer valid for $m=2$. In this case, the linear system (\ref{eq:linear_system}) reduces to a single equation that can be solved to yield:
\begin{align*}
G_{1}(s) &= (B_{1}(\ell_{1},s) - A_{2}(\ell_{1},s))^{-1}\bigl[P_{2}(\ell_{1},s) - P_{1}(\ell_{1},s) - A_{1}(\ell_{1},s)G_{0}(s) + B_{2}(\ell_{1},s)G_{L}(s)\bigr],
\end{align*}
which allows $C_{1}(x,s)$ and $C_{2}(x,s)$ to be evaluated at any $x$ and $s$.

\subsection{Numerical inversion of Laplace transform}
\label{sec:solution3}
To convert the Laplace domain expressions (\ref{eq:C1})--(\ref{eq:Cm}) back to the time domain and thereby obtain the concentration $c_{i}(x,t)$ ($i = 1,\hdots,m$), we need to compute the inverse Laplace transform:
\begin{align*}
c_{i}(x,t) &= \mathcal{L}^{-1}\left\{C_{i}(x,s)\right\} = \frac{1}{2\pi i}\int_{\Gamma} e^{st}C_{i}(x,s)\,\text{d}s,
\end{align*}
where $\Gamma$ is a Hankel contour that begins at $-\infty-0i$, winds around the origin and terminates at $-\infty+0i$ {\citep{trefethen_2006}}. Introducing the change of variable $z = st$, using a rational approximation to $e^{z}$ and applying residue calculus (see {\cite{trefethen_2006}} for full details) yields the numerical inversion formula:
\begin{align}
\label{eq:trefethen}
c_{i}(x,t) &= \mathcal{L}^{-1}\left\{C_{i}(x,s)\right\} \approx - \frac{2}{t}\Re\Biggl\{\sum_{k\in\mathbb{O}_{N}} w_{k}C_{i}\left(x,s_{k}\right)\Biggr\},\hspace{-1cm}
\end{align}
where $N$ is even, $\mathbb{O}_{N}$ is the set of positive odd integers less than $N$, $s_{k} = z_{k}/t$ and $w_{k},z_{k}\in\mathbb{C}$ are the residues and poles of the best $(N,N)$ rational approximation to $e^{z}$ on the negative real line. Both $w_{k}$ and $z_{k}$ are constants, which are independent of $x$ and $t$ and computed using a supplied MATLAB function \citep[Fig~4.1]{trefethen_2006}.

\subsection{Treatment of step function boundary conditions}
Suppose the inlet concentration $c_{0}(t)$ in either the concentration-type (\ref{eq:inlet_concentration_type}) or flux type (\ref{eq:inlet_flux_type}) boundary condition, is a Heaviside step function of duration $t_{0}>0$:
\begin{gather*}
c_{0}(t) = c_{0}H(t_{0}-t) = \begin{cases} c_{0}, & 0 < t < t_{0},\\ 0, & t > t_{0},\end{cases}
\end{gather*}
where $c_{0}$ is a constant. In solute transport problems, such boundary conditions are commonly \citep{van_genuchten_1982,leij_1991,goltz_2017} paired with a zero concentration gradient at the outlet (\ref{eq:outlet_zero_concentration_gradient}) and lead to $G_{0}(s) = \exp(-t_{0}s)/s$ and $G_{0}(s) = v_{1}\exp(-t_{0}s)/s$ in Eq (\ref{eq:C1}) for the concentration-type and flux-type boundary condition, respectively. Such exponential functions are well known to cause numerical problems in algorithms for inverting Laplace transforms \citep{kuhlman_2013}. The approximation (\ref{eq:trefethen}) indeed suffers from this issue as evaluating $\exp(-t_{0}s)$ at $s = s_{k} = z_{k}/t$ for poles $z_{k}$ with negative real part leads to floating point overflow for small $t$. To overcome this problem, we use superposition of solutions. Consider, for example, the multilayer transport model (\ref{eq:pde})--(\ref{eq:dispersive_flux_continuity}) subject to the boundary conditions (\ref{eq:inlet_concentration_type}) or (\ref{eq:inlet_flux_type}), and (\ref{eq:outlet_zero_concentration_gradient}). The solution to this problem can be expressed as
\begin{align*}
c_{i}(x,t) = \begin{cases} \widetilde{c}_{i}(x,t), & 0 < t < t_{0},\\
\widetilde{c}_{i}(x,t) - \widehat{c}_{i}(x,t-t_{0}), & t > t_{0},\end{cases}
\end{align*}
where $\widetilde{c}_{i}(x,t)$ is the solution of the multilayer transport model (\ref{eq:pde})--(\ref{eq:general_bcR}) with $g_{0}(t) = c_{0}$ and $\widehat{c}_{i}(x,t)$ is the solution of the multilayer transport model (\ref{eq:pde})--(\ref{eq:general_bcR}) with $g_{0}(t) = c_{0}$, $f_{i} = 0$ and $\gamma_{i} = 0$. Both $\widetilde{c}_{i}(x,t)$ and  $\widehat{c}_{i}(x,t)$ are obtained using the approach outlined in subsections \ref{sec:solution1}--\ref{sec:solution3}. 

\section{Results}
\label{sec:results}
We now demonstrate application of our semi-analytical Laplace-transform solution and verify that it produces the correct results using a selection of test cases. The transport parameters, initial conditions and boundary conditions for each problem are provided in Tables \ref{tab:transport_parameters} and \ref{tab:boundary_conditions}. A MATLAB code implementing our semi-analytical solution and producing the results in this section is available for download from \href{https://github.com/elliotcarr/Carr2020a}{github.com/elliotcarr/Carr2020a}.

\begin{table*}
\centering\small
\caption{Transport parameters, geometry and initial conditions for the various test cases, where $i$ is the layer index and $\ell_{i}$ locates the end of layer $i$ [$\textrm{cm}$]. In the $i$th layer: $R_{i}$ is the retardation factor $[-]$, $D_{i}$ is the dispersion coefficient [$\textrm{cm}^{2}\text{day}^{-1}$], $v_{i}$ is the pore-water velocity [$\textrm{cm}\,\text{day}^{-1}$], $\mu_{i}$ is the rate constant for first-order decay [$\textrm{day}^{-1}$], $\gamma_{i}$ is the rate constant for zero-order production [$\text{kg}\,\text{cm}^{-3}\,\text{day}^{-1}$], $\theta_{i}$ is the volumetric water content $[-]$ and $f_{i}$ is the initial constant concentration [$\textrm{kg}\,\text{cm}^{-3}$]. Dashed horizontal lines for case 12 indicate the artificial layer created to accommodate the injected contaminant initial condition (see section \ref{sec:results}).}
\begin{tabular*}{0.6\textwidth}{@{\extracolsep{\fill}}lrrrrrrrrr}
\hline
Case & $i$ & $\ell_{i}$ & $R_{i}$ & $D_{i}$ & $v_{i}$ & $\mu_{i}$ & $\gamma_{i}$ & $\theta_{i}$ & $f_{i}$\\
\hline
1--4 & 1 & 10 & 1 & 50 & 75 & 2 & 1 & 0.4 & 0\\
  & 2 & 30 & 1 & 50 & 75 & 2 & 1 & 0.4 & 0\\
\hline
5 & 1 & 10 & 1 & 50 & 25 & 0 & 0 & 0.4 & 0\\
  & 2 & 30 & 1 & 20 & 40 & 0 & 0 & 0.25 & 0\\
\hline  
6 & 1 & 10 & 1 & 50 & 25 & 0 & 0 & 0.4 & 0\\
  & 2 & 30 & 1 & 20 & 40 & 0 & 0 & 0.25 & 0\\
\hline  
7 & 1 & 10 & 1 & 50 & 25 & 0 & 0 & 0.4 & 0\\
  & 2 & 30 & 1 & 20 & 40 & 0 & 0 & 0.25 & 0\\ 
\hline
8 & 1 & 10 & 3 & 50 & 25 & 3 & 0 & 0.4 & 0\\
  & 2 & 20 & 2 & 20 & 40 & 4 & 0 & 0.25 & 0\\         
\hline
9-10, 11 & 1 & 10 & 4.25 & 7 & 10 & 0 & 0 & 0.4 & 0\\ 
  & 2 & 12 & 14 & 18 & 8 & 0 & 0 & 0.5 & 0\\
  & 3 & 20 & 4.25 & 7 & 10 & 0 & 0 & 0.4 & 0\\ 
  & 4 & 22 & 14 & 18 & 8 & 0 & 0 & 0.5 & 0\\
  & 5 & 30 & 4.25 & 7 & 10 & 0 & 0 & 0.4 & 0\\
\hline
12 & 1 & 10 & 4.25 & 7 & 10 & 0 & 0 & 0.4 & 0\\ 
  & 2 & 12 & 14 & 18 & 8 & 0 & 0 & 0.5 & 0\\
  \cdashline{2-10}
  & 3 & 14 & 4.25 & 7 & 10 & 0 & 0 & 0.4 & 0\\ 
 & 4 & 18 & 4.25 & 7 & 10 & 0 & 0 & 0.4 & $c_{0}$\\ 
  & 5 & 20 & 4.25 & 7 & 10 & 0 & 0 & 0.4 & 0\\
  \cdashline{2-10} 
  & 6 & 22 & 14 & 18 & 8 & 0 & 0 & 0.5 & 0\\
  & 7 & 30 & 4.25 & 7 & 10 & 0 & 0 & 0.4 & 0\\ 
\hline
13 & 1 & 10 & 4.25 & 7 & 10 & 3 & 2 & 0.4 & 0\\ 
  & 2 & 12 & 14 & 18 & 8 & 2 & 4 & 0.5 & 0\\
  & 3 & 20 & 4.25 & 7 & 10 & 3 & 2 & 0.4 & 0\\ 
  & 4 & 22 & 14 & 18 & 8 & 2 & 4 & 0.5 & $c_{0}$\\
  & 5 & 30 & 4.25 & 7 & 10 & 3 & 2 & 0.4 & 0\\
\hline
\end{tabular*}
\label{tab:transport_parameters}
\end{table*}

\begin{table*}
\centering\small
\caption{Choice of inlet boundary condition (\ref{eq:general_bcL}) for the various test cases. For all test cases, a zero concentration gradient is assumed at the outlet with parameters $a_{L} = 0$, $b_{L} = 1$ and $g_{L}(t) = 0$ specified in the outlet boundary condition (\ref{eq:general_bcR}).}
\begin{tabular*}{0.5\textwidth}{@{\extracolsep{\fill}}lrrr}
\hline
Case & $a_{0}$ & $b_{0}$ & $g_{0}(t)$\\
\hline
1, 5--9 & $v_{1}$ & $D_{1}$ & $v_{1}c_{0}$\\
2, 10, 13 & $v_{1}$ & $D_{1}$ & $v_{1}c_{0}H(t_{0}-t)$\\
3 & 1 & 0 & $c_{0}$\\
4 & 1 & 0 & $c_{0}H(t_{0}-t)$\\
11 & 1 & 0 & $c_{0}\alpha t e^{-\beta t}$\\
12 & 0 & 1 & 0\\
\hline
\end{tabular*}
\label{tab:boundary_conditions}
\end{table*}

\begin{table*}
\centering\small
\caption{Comparison between our semi-analytical Laplace transform solution (see section \ref{sec:solution}) and the analytical solutions catalogued in \cite{van_genuchten_1982} for the homogeneous medium test cases (1--4). The tabulated values are the maximum absolute difference between the values of the relative concentration ($c(x,t)/c_{0}$) computed using the two approaches over the discrete range $x = 0,2,\hdots,20\,\mathrm{cm}$, with $N = 14$ poles/residues used for numerically inverting the Laplace transform (\ref{eq:trefethen}). As recommended by \citet {van_genuchten_1982}, if $v_{1}L/D_{1} > \min(5+40v_{1}t/(R_{1}L),5+40v_{1}(t-t_{0})/(R_{1}L),100)$, the analytical solution is calculated using the approximate solutions catalogued \citep{van_genuchten_1982}. Otherwise the eigenfunction expansion solutions are used with $1000$ eigenvalues/terms taken in the series expansions (see \citet{van_genuchten_1982} and our code for full details).}
\begin{tabular*}{0.95\textwidth}{@{\extracolsep{\fill}}cllllll}
\hline
Case & $t = 10^{-3}$ & $t = 0.1$ & $t = 0.6$ & $t = 1$ & $t = 2$ & $t = 4$\\
\hline
1 & \num{4.11e-14} & \num{5.53e-10} & \num{5.84e-08} & \num{9.38e-09} & \num{4.75e-09} & \num{6.10e-10}\\
2 & \num{4.11e-14} & \num{5.53e-10} & \num{1.89e-08} & \num{7.10e-08} & \num{1.74e-08} & \num{5.90e-10}\\
3 & \num{1.98e-14} & \num{7.61e-10} & \num{1.93e-08} & \num{1.53e-08} & \num{1.20e-09} & \num{3.34e-10}\\
4 & \num{1.98e-14} & \num{7.61e-10} & \num{1.93e-08} & \num{1.57e-09} & \num{1.15e-08} & \num{8.10e-10}\\
\hline
\end{tabular*}
\label{tab:Case1}
\end{table*}

\subsection{One and two layer test cases}
Cases 1--4 consider advection-dispersion-reaction in a homogeneous (single-layer) medium of length $30\,\text{cm}$. The concentration is initially zero everywhere and a zero concentration gradient is applied at the outlet (\ref{eq:outlet_zero_concentration_gradient}). Four different boundary conditions are considered at the inlet: 
\begin{itemize}
\item Case 1: flux-type boundary condition (\ref{eq:inlet_flux_type}) with constant inlet concentration $c_{0}(t) = c_{0}$; 
\item Case 2: flux-type boundary condition (\ref{eq:inlet_flux_type}) with step inlet concentration $c_{0}(t) = c_{0}H(t_{0}-t)$ and pulse duration $t_{0}=0.5\,\textrm{days}$; 
\item Case 3: concentration-type boundary condition  (\ref{eq:inlet_concentration_type}) with constant inlet concentration $c_{0}(t) = c_{0}$;
\item Case 4: concentration-type boundary condition  (\ref{eq:inlet_concentration_type}) with step inlet concentration $c_{0}(t) = c_{0}H(t_{0}-t)$ and pulse duration $t_{0}=0.5\,\textrm{days}$.
\end{itemize}
To solve these single-layer problems using our semi-analytical method, we choose $m=2$ layers and set the transport parameters equal in both layers (see Table \ref{tab:transport_parameters}). In Table \ref{tab:Case1}, the relative concentration distributions ($c(x,t)/c_{0}$) obtained are compared to corresponding distributions obtained using analytical solutions given in the literature \cite[sections C7 and C8]{van_genuchten_1982}, which are valid for homogeneous (single-layer) media only. For all four problems, the maximum absolute difference between the relative concentration distributions obtained using both solution methods is tabulated in Table \ref{tab:Case1} at $t = 10^{-3},0.1,0.6,1,2,4\,\textrm{days}$. These results demonstrate that both solutions are in excellent agreement for each choice of inlet condition and verify that our semi-analytical solution produces the correct results for single-layer media.

\begin{table*}[p]
\small\centering
\caption{Relative concentration values ($c(x,t)/c_{0}$) for test cases 5--7 computed using our semi-analytical Laplace transform (SALT) solution (see section \ref{sec:solution}) with $N = 14$ poles/residues used in the numerical inversion of the Laplace transform (\ref{eq:trefethen}). Results are compared with previously published values given by \citet*{guerrero_2013} (GPS), \citet*{liu_1998} (LBE) and \citet*{leij_1995} (LVG). Shaded cells highlight discrepancies between the different methods by indicating when one of the four computed values differs from the others.}
\setlength{\tabcolsep}{0.1cm}
\renewcommand*{\arraystretch}{1.0}
\begin{tabular*}{1.0\textwidth}{|c|r|rrrr|rrrr|rrrr|rrrr|rrrr|}
\cline{1-18}
{\footnotesize Case}\hspace*{-0.03cm} & & \multicolumn{4}{l|}{$t = 0.2$} & \multicolumn{4}{l|}{$t = 0.4$} & \multicolumn{4}{l|}{$t = 0.6$} & \multicolumn{4}{l|}{$t = 0.8$}\\
\cline{3-18}
& $x$ & {\footnotesize SALT} & GPS & LBE & LVG & {\footnotesize SALT} & GPS & LBE & LVG & {\footnotesize SALT} & GPS & LBE & LVG & {\footnotesize SALT} & GPS & LBE & LVG\\
\cline{1-18}
5 & 0 & 0.884 & 0.884 & 0.884 & 0.884 & 0.963 & 0.963 & 0.963 & 0.963 & 0.987 & 0.987 & 0.987 & 0.987 & 0.995 & 0.995 & 0.995 & 0.995\\
 & 2 & 0.742 & 0.742 & 0.742 & 0.742 & 0.915 & 0.915 & 0.915 & 0.915 & 0.969 & 0.969 & 0.969 & 0.969 & 0.988 & 0.988 & 0.988 & 0.988\\
 & 4 & 0.561 & 0.561 & 0.561 & 0.561 & 0.841 & 0.841 & 0.841 & 0.841 & 0.940 & 0.940 & 0.940 & 0.940 & 0.977 & 0.977 & 0.977 & 0.977\\
 & 6 & 0.375 & 0.375 & \cellcolor{tableshade}{0.374} & 0.375 & 0.746 & 0.746 & 0.746 & 0.746 & 0.901 & 0.901 & 0.901 & 0.901 & 0.962 & 0.962 & 0.962 & 0.962\\
 & 8 & 0.222 & 0.222 & 0.222 & 0.222 & 0.645 & 0.645 & 0.645 & 0.645 & 0.858 & 0.858 & 0.858 & 0.858 & 0.945 & 0.945 & 0.945 & 0.945\\
 & 10 & 0.142 & 0.142 & 0.142 & 0.142 & 0.579 & 0.579 & 0.579 & 0.579 & 0.829 & 0.829 & 0.829 & 0.829 & 0.933 & 0.933 & 0.933 & 0.933\\
 & 12 & 0.063 & 0.063 & 0.063 & 0.063 & 0.480 & 0.480 & 0.480 & 0.480 & 0.781 & 0.781 & 0.781 & 0.781 & 0.914 & 0.914 & 0.914 & 0.914\\
 & 14 & 0.021 & 0.021 & 0.021 & 0.021 & 0.372 & 0.372 & 0.372 & 0.372 & 0.722 & 0.722 & 0.722 & 0.722 & 0.889 & 0.889 & 0.889 & 0.889\\
 & 16 & 0.005 & 0.005 & 0.005 & 0.005 & 0.264 & 0.264 & \cellcolor{tableshade}{0.265} & 0.264 & 0.651 & 0.651 & 0.651 & 0.651 & 0.858 & 0.858 & 0.858 & 0.858\\
 & 18 & 0.001 & 0.001 & 0.001 & 0.001 & 0.168 & 0.168 & \cellcolor{tableshade}{0.169} & 0.168 & 0.567 & 0.567 & 0.567 & 0.567 & 0.819 & 0.819 & 0.819 & 0.819\\
 & 20 & 0.000 & 0.000 & 0.000 & 0.000 & 0.094 & 0.094 & 0.094 & 0.094 & 0.473 & 0.473 & 0.473 & 0.473 & 0.770 & 0.770 & 0.770 & 0.770\\
\cline{1-18}
6 & 0 & 0.978 & 0.978 & \cellcolor{tableshade}{0.977} & 0.978 & 0.998 & 0.998 & 0.998 & 0.998 & 1.000 & 1.000 & 1.000 & 1.000 & 1.000 & 1.000 & 1.000 & 1.000\\
 & 2 & 0.868 & 0.868 & \cellcolor{tableshade}{0.867} & 0.868 & 0.984 & 0.984 & 0.984 & 0.984 & 0.998 & 0.998 & 0.998 & 0.998 & 1.000 & 1.000 & 1.000 & 1.000\\
 & 4 & 0.634 & 0.634 & \cellcolor{tableshade}{0.633} & 0.634 & 0.942 & 0.942 & 0.942 & 0.942 & 0.991 & 0.991 & 0.991 & 0.991 & 0.999 & 0.999 & 0.999 & 0.999\\
 & 6 & 0.345 & 0.345 & 0.345 & 0.345 & 0.849 & 0.849 & 0.849 & 0.849 & 0.972 & 0.972 & 0.972 & 0.972 & 0.995 & 0.995 & 0.995 & 0.995\\
 & 8 & 0.131 & 0.131 & 0.131 & 0.131 & 0.693 & 0.693 & 0.693 & 0.693 & 0.930 & 0.930 & \cellcolor{tableshade}{0.929} & 0.930 & 0.986 & 0.986 & 0.986 & 0.986\\
 & 10 & 0.033 & 0.033 & 0.033 & 0.033 & 0.496 & 0.496 & 0.496 & 0.496 & 0.853 & 0.853 & 0.853 & 0.853 & 0.966 & 0.966 & 0.966 & 0.966\\
 & 12 & 0.011 & 0.011 & 0.011 & 0.011 & 0.370 & 0.370 & 0.370 & 0.370 & 0.784 & 0.784 & \cellcolor{tableshade}{0.783} & 0.784 & 0.944 & 0.944 & 0.944 & 0.944\\
 & 14 & 0.003 & 0.003 & 0.003 & 0.003 & 0.257 & 0.257 & 0.257 & 0.257 & 0.699 & 0.699 & \cellcolor{tableshade}{0.698} & 0.699 & 0.913 & 0.913 & 0.913 & 0.913\\
 & 16 & 0.001 & 0.001 & 0.001 & 0.001 & 0.166 & 0.166 & 0.166 & 0.166 & 0.601 & 0.601 & 0.601 & 0.601 & 0.871 & 0.871 & 0.871 & 0.871\\
 & 18 & 0.000 & 0.000 & 0.000 & 0.000 & 0.098 & 0.098 & \cellcolor{tableshade}{0.099} & 0.098 & 0.498 & 0.498 & 0.498 & 0.498 & 0.817 & 0.817 & 0.817 & 0.817\\
 & 20 & 0.000 & 0.000 & 0.000 & 0.000 & 0.054 & 0.054 & 0.054 & 0.054 & 0.395 & 0.395 & 0.395 & 0.395 & 0.751 & 0.751 & \cellcolor{tableshade}{0.750} & 0.751\\
\cline{1-18}
7 & 0 & 0.999 & 0.999 & 0.999 & 0.999 & 1.000 & 1.000 & 1.000 & 1.000 & 1.000 & 1.000 & 1.000 & 1.000 & 1.000 & 1.000 & 1.000 & 1.000\\
 & 2 & 0.988 & 0.988 & \cellcolor{tableshade}{0.987} & 0.988 & 1.000 & 1.000 & 1.000 & 1.000 & 1.000 & 1.000 & 1.000 & 1.000 & 1.000 & 1.000 & 1.000 & 1.000\\
 & 4 & 0.928 & 0.928 & 0.928 & 0.928 & 0.999 & 0.999 & 0.999 & 0.999 & 1.000 & 1.000 & 1.000 & 1.000 & 1.000 & 1.000 & 1.000 & 1.000\\
 & 6 & 0.764 & 0.764 & \cellcolor{tableshade}{0.763} & 0.764 & 0.995 & 0.995 & 0.995 & 0.995 & 1.000 & 1.000 & 1.000 & 1.000 & 1.000 & 1.000 & 1.000 & 1.000\\
 & 8 & 0.496 & 0.496 & \cellcolor{tableshade}{0.495} & 0.496 & 0.976 & 0.976 & 0.976 & 0.976 & 0.998 & 0.998 & 0.998 & 0.998 & 0.999 & 0.999 & 0.999 & 0.999\\
 & 10 & 0.152 & 0.152 & 0.152 & 0.152 & 0.780 & 0.780 & \cellcolor{tableshade}{0.779} & 0.780 & 0.940 & 0.940 & \cellcolor{tableshade}{0.939} & 0.940 & 0.979 & 0.979 & \cellcolor{tableshade}{0.978} & 0.979\\
 & 12 & 0.049 & 0.049 & \cellcolor{tableshade}{0.050} & 0.049 & 0.600 & 0.600 & 0.600 & 0.600 & 0.870 & 0.870 & 0.870 & 0.870 & 0.952 & 0.952 & 0.952 & 0.952\\
 & 14 & 0.013 & 0.013 & 0.013 & 0.013 & 0.418 & 0.418 & 0.418 & \cellcolor{tableshade}{0.417} & 0.773 & 0.773 & 0.773 & 0.773 & 0.911 & 0.911 & \cellcolor{tableshade}{0.910} & 0.911\\
 & 16 & 0.003 & 0.003 & 0.003 & 0.003 & 0.262 & 0.262 & 0.262 & 0.262 & 0.653 & 0.653 & 0.653 & 0.653 & 0.851 & 0.851 & 0.851 & 0.851\\
 & 18 & 0.000 & 0.000 & 0.000 & 0.000 & 0.148 & 0.148 & 0.148 & 0.148 & 0.522 & 0.522 & 0.522 & 0.522 & 0.774 & 0.774 & 0.774 & 0.774\\
 & 20 & 0.000 & 0.000 & 0.000 & 0.000 & 0.075 & 0.075 & 0.075 & 0.075 & 0.393 & 0.393 & 0.393 & 0.393 & 0.681 & 0.681 & 0.681 & 0.681\\
\cline{1-18}
\end{tabular*}
\label{tab:Case2}
\end{table*}

\begin{figure*}[p]
\centering
\includegraphics[width=\figw]{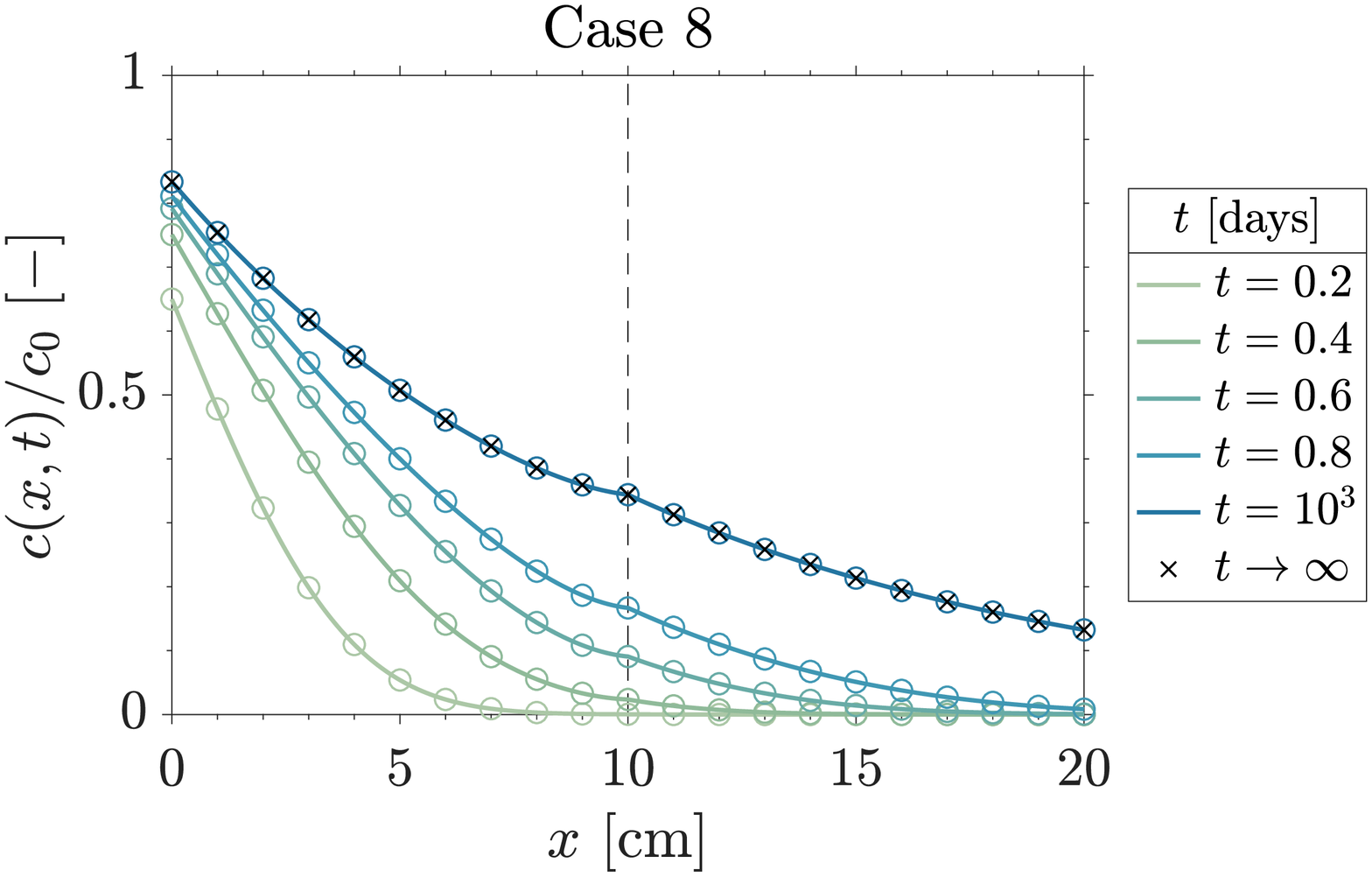}
\caption{Relative concentration spatial profiles ($c(x,t)/c_{0}$) over time for test case 8. Continuous lines indicate our semi-analytical Laplace transform solution (see section \ref{sec:solution}) with $N = 14$ poles/residues used in the numerical inversion of the Laplace transform (\ref{eq:trefethen}). Circle markers indicate the numerical solution obtained using the finite volume scheme outlined in Appendix \ref{app:fvm} computed with $n = 601$ nodes but plotted with markers at $x = 0,1,\hdots,20\,\textrm{cm}$ only. Crosses indicate the exact steady-state solution obtained by solving the steady-state analogue of the multilayer transport model (\ref{eq:pde})--(\ref{eq:general_bcR}) analytically. The vertical dashed line indicates the location of the interface at $x = 10\,\textrm{cm}$. The plot is truncated at $x = 20\,\textrm{cm}$ rather than shown for the full length of the medium ($x = 30\,\mathrm{cm}$) inline with \citet{guerrero_2013}.}
\label{fig:plots1}
\end{figure*}

\begin{figure*}[!b]
\centering
\includegraphics[width=\figw]{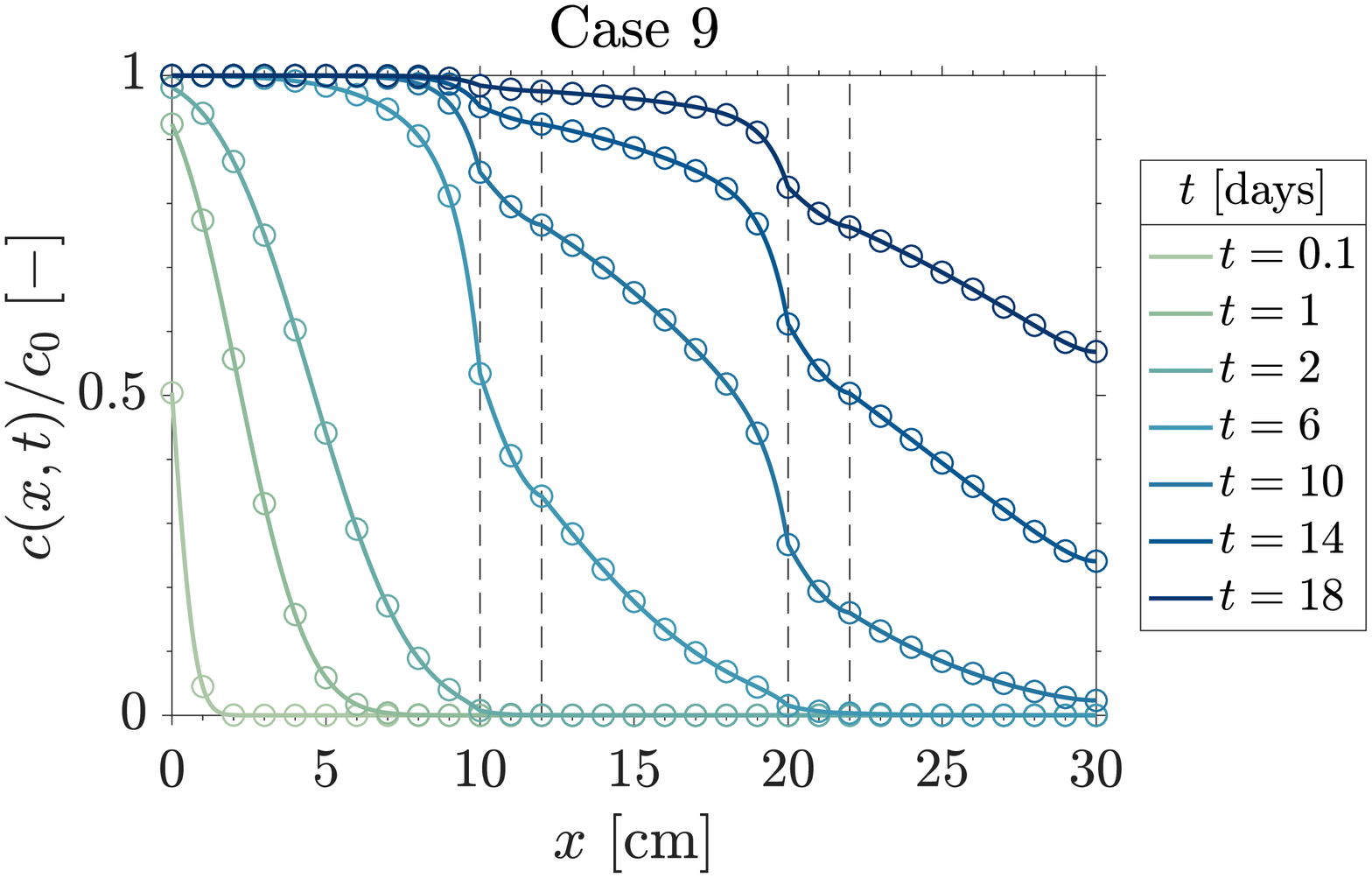}\hspace{0.01\textwidth}
\includegraphics[width=\figw]{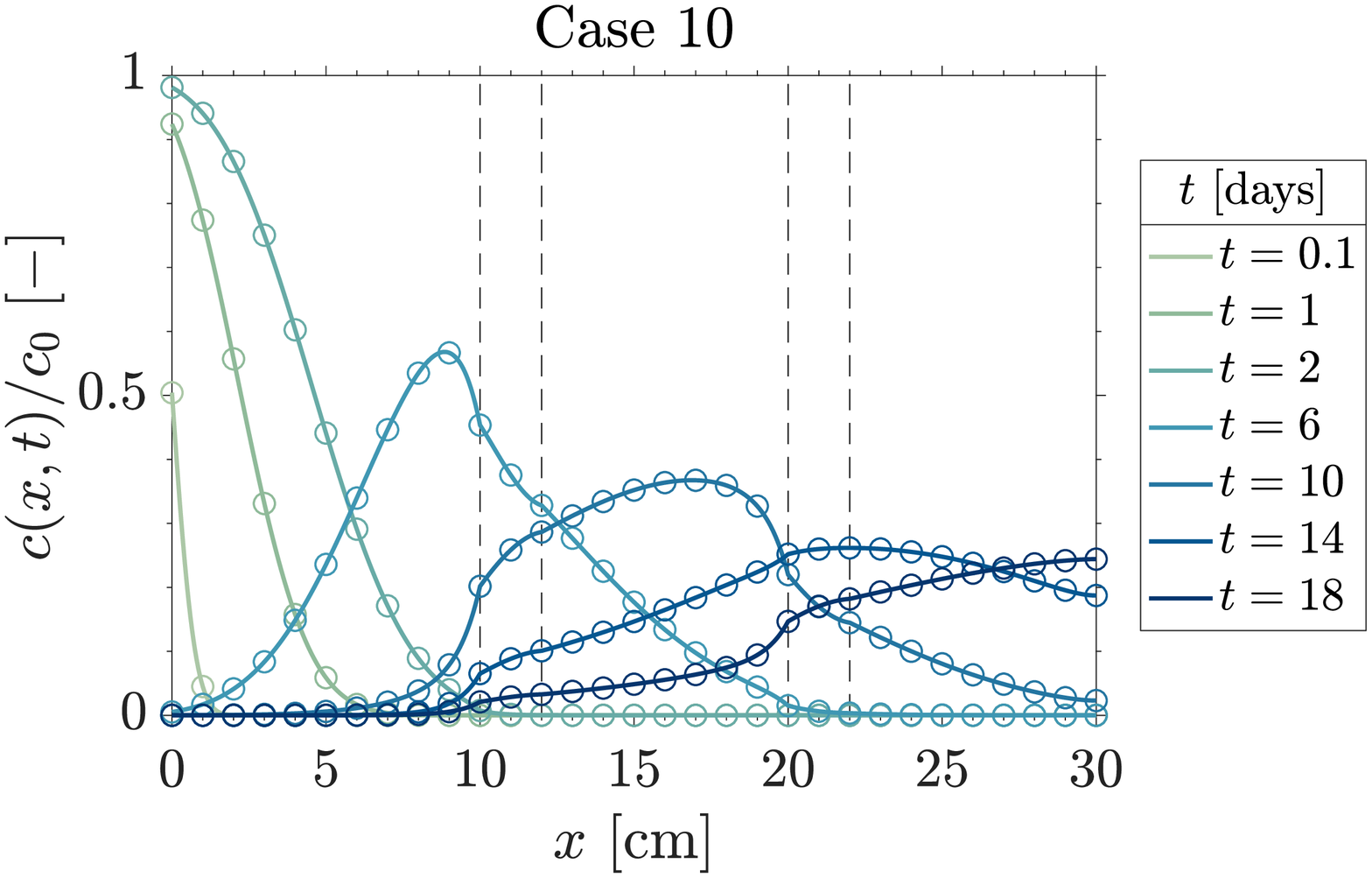}\\[0.5cm]
\includegraphics[width=\figw]{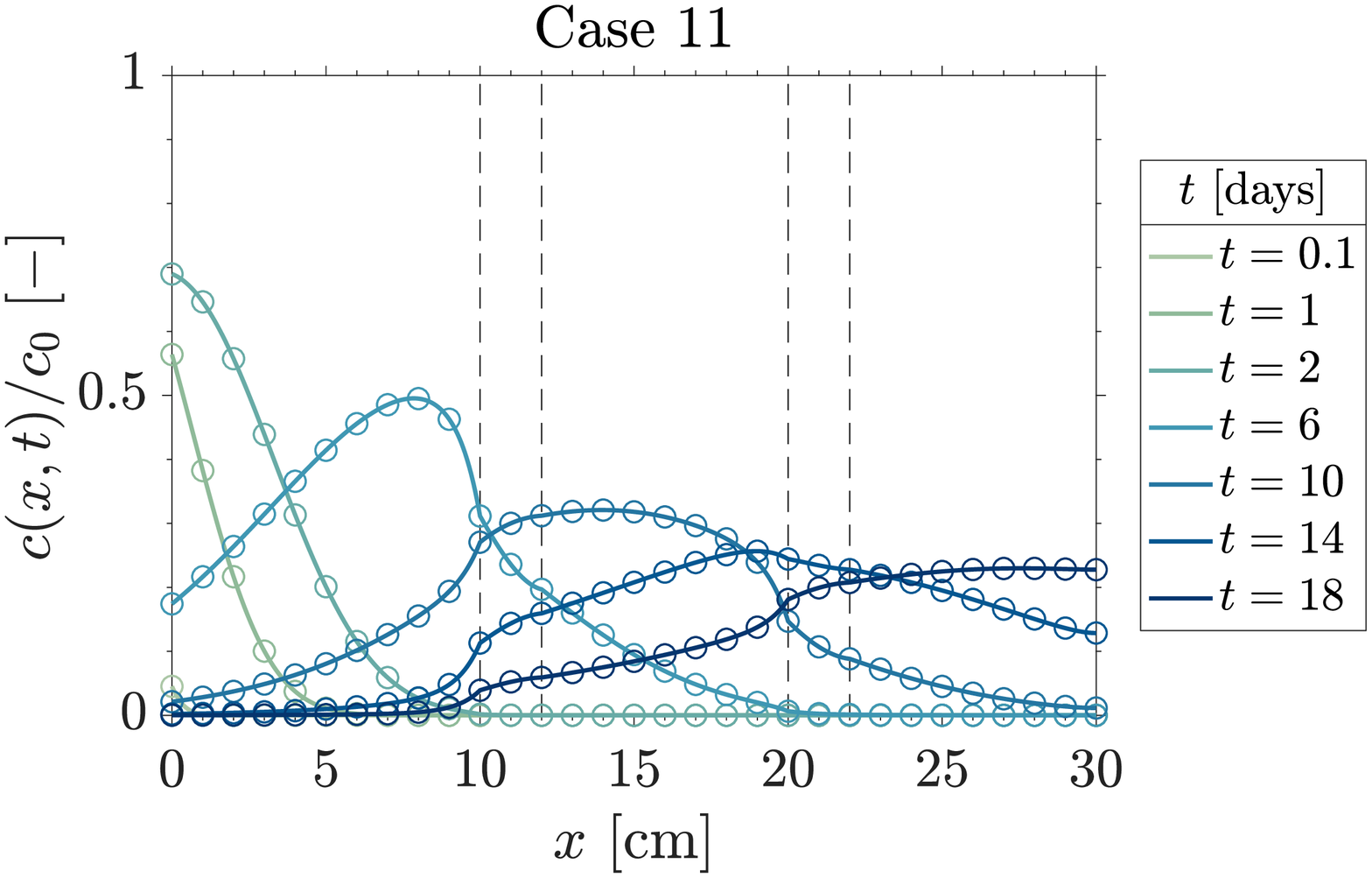}\hspace{0.01\textwidth}
\includegraphics[width=\figw]{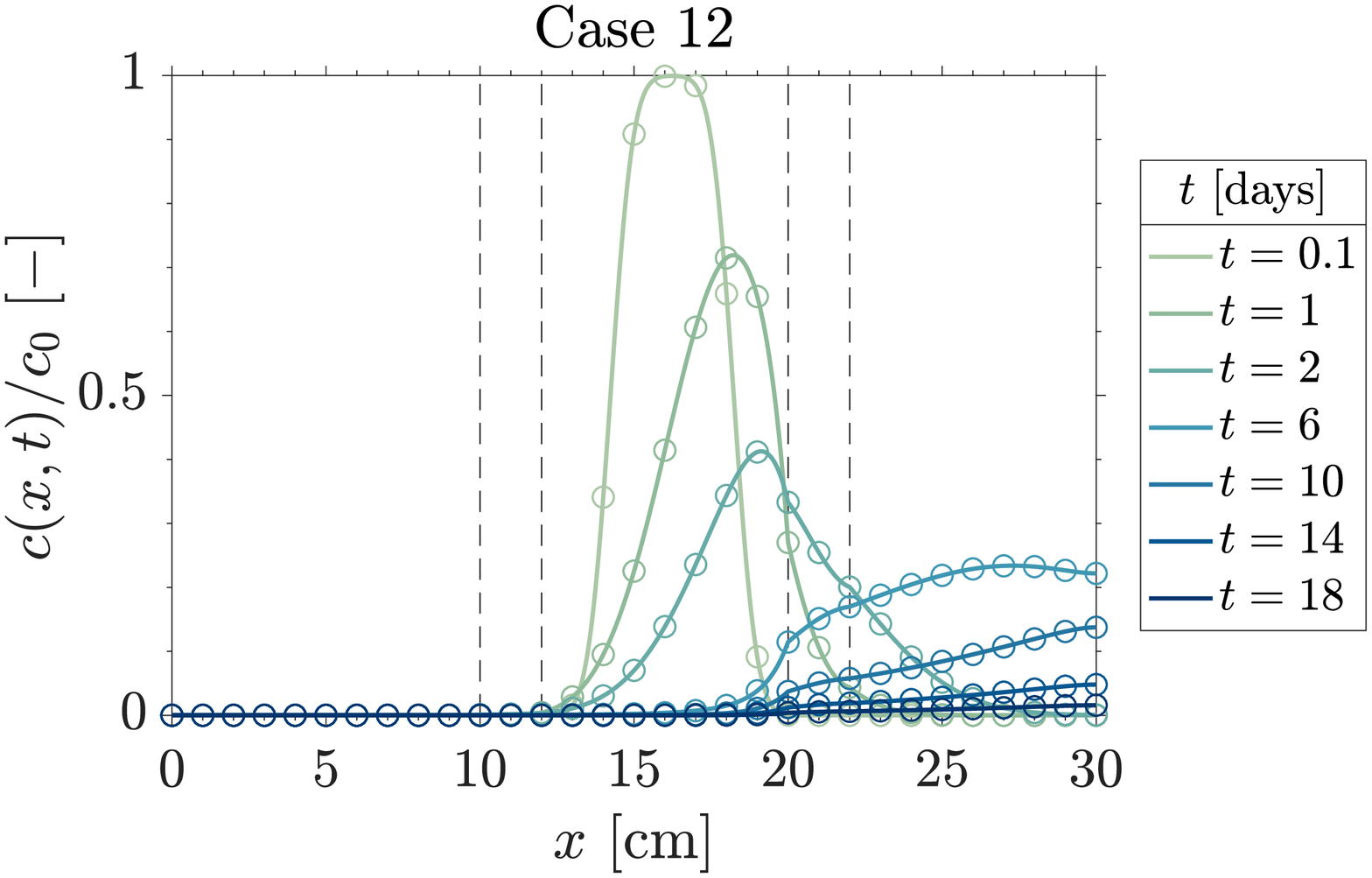}\\[0.5cm]
\includegraphics[width=\figw]{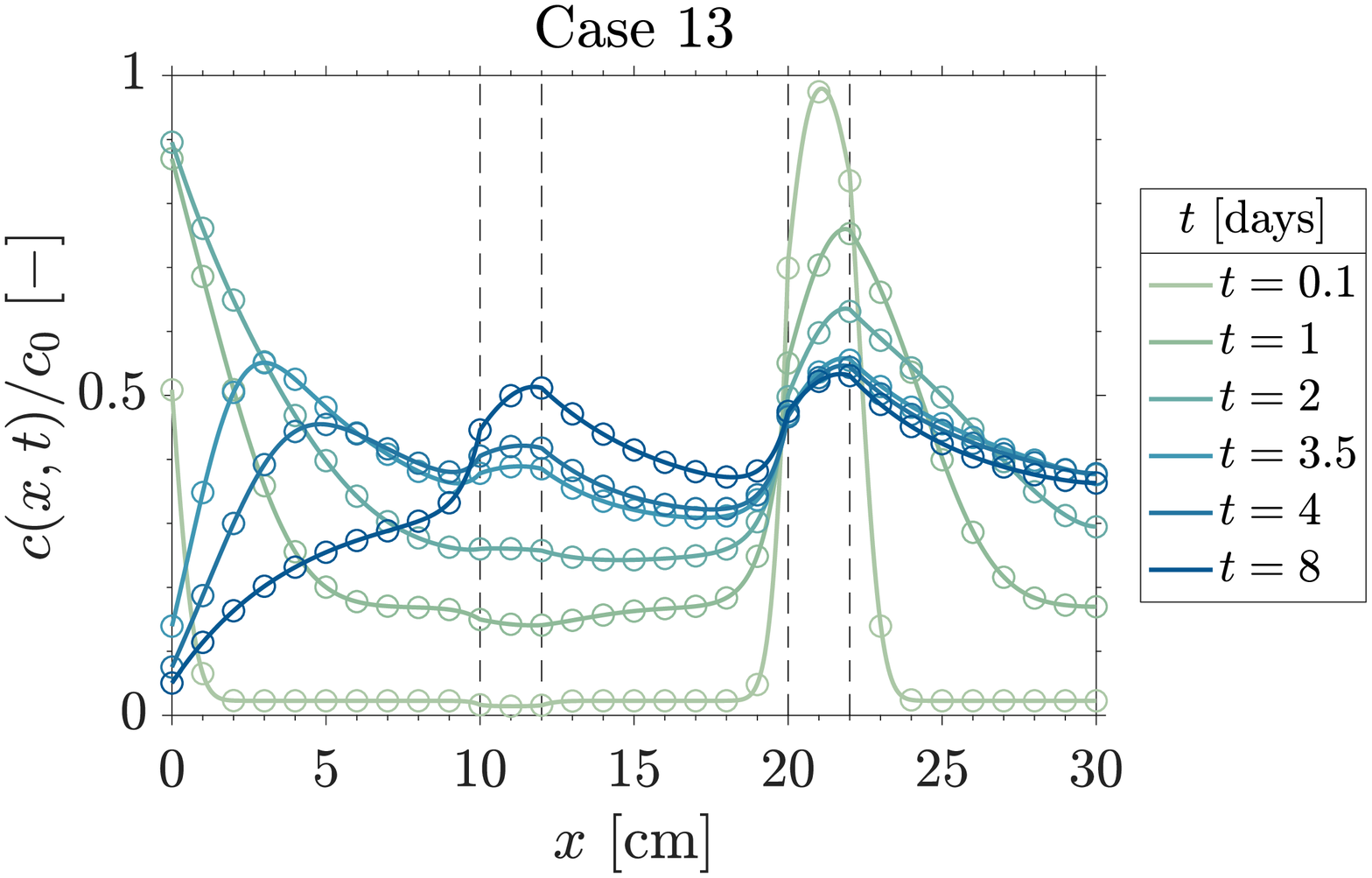}
\caption{Relative concentration spatial profiles ($c(x,t)/c_{0}$) over time for test cases 9--13. Continuous lines indicate our semi-analytical Laplace transform solution (see section \ref{sec:solution}) with $N = 14$ poles/residues used in the numerical inversion of the Laplace transform (\ref{eq:trefethen}). Circle markers indicate the numerical solution obtained using the finite volume scheme outlined in Appendix \ref{app:fvm} computed with $n = 601$ nodes but plotted with markers at $x = 0,1,\hdots,30\,\textrm{cm}$ only. Vertical dashed lines indicate the location of the layer interfaces at $x = 10,12,20,22\,\textrm{cm}$.}
\label{fig:plots2}
\end{figure*}

Next, we present results for three two-layer test cases that have frequently appeared in the literature \citep{leij_1995,liu_1998,guerrero_2013}. The three test cases, labelled cases 5--7, consider advection-dispersion (without decay or production) in a medium of length $30\,\textrm{cm}$ with first and second layers of length $10\,\textrm{cm}$ and $20\,\textrm{cm}$, respectively. A constant flux-type boundary condition is applied at the inlet (Eq (\ref{eq:inlet_flux_type}) with $c_{0}(t) = c_{0}$) and a zero concentration gradient is applied at the outlet (\ref{eq:outlet_zero_concentration_gradient}). Initially, the concentration is assumed to be zero in both layers. Different transport parameter combinations are applied for the three test cases (see Table \ref{tab:transport_parameters}). In Table \ref{tab:Case2}, for all three test cases, we report the relative concentration values obtained from our semi-analytical solution at several equidistant points in space and time. These results are compared to those previously reported by \citet{guerrero_2013}, \citet{liu_1998} and \citet{leij_1995} with all numerical values in Table \ref{tab:Case2} displayed to three decimal places to be consistent with the numerical values reported in those papers. For all three test cases, our results agree precisely with those of \citet{guerrero_2013} to the precision reported. When compared to \citet{liu_1998} and \citet{leij_1995}, minor differences of $\pm 0.001$ are evident as highlighted in Table \ref{tab:Case2}. These discrepancies are likely to be explained by \citet{leij_1995} considering a semi-infinite second layer and using a different method for numerically inverting the Laplace transform and \citet{liu_1998} considering a finite second layer of undisclosed length.

Lastly, we present results for a final two-layer test case previously considered by \citet{guerrero_2013}. This problem, labelled case 8, is the same as case 6 with the exception that first-order decay is present in both layers (see Table \ref{tab:transport_parameters}). In Figure \ref{fig:plots1}, we plot the relative concentration distributions obtained from our semi-analytical solution at $t = 0.2,0.4,0.6,0.8\,\textrm{days}$ (as in \citet{guerrero_2013}) and $t = 10^{3}$ when the solution is visible indistinguishable from its steady state. At short time, our results match well with \citet{guerrero_2013}'s results but differ for long times and at steady-state. However, agreement with the numerical solution (outlined in Appendix \ref{app:fvm}) as well as the exact steady-state solution obtained by solving the steady-state analogue of the multilayer transport model (\ref{eq:pde})--(\ref{eq:general_bcR}) (see Figure \ref{fig:plots1}) supports that our approach produces the correct results.

\subsection{Multiple layer test cases}
We now apply our semi-analytical solution to some test cases involving five or more layers. Our solutions are compared to previously reported results from the literature and/or numerical solutions obtained using a finite volume method, briefly discussed in Appendix \ref{app:fvm}.

Firstly, we consider the five-layer test case previously considered by \citet{liu_1998} and \citet{guerrero_2013}. This test case involves advection-dispersion (without decay or production) in a medium of length $30\,\mathrm{cm}$ consisting of alternating layers of sand--clay--sand--clay--sand (cases 9--10 in Tables \ref{tab:transport_parameters} and \ref{tab:boundary_conditions}). The initial concentration is zero in each layer and a zero concentration gradient is applied at the outlet. Two different boundary conditions are considered at the inlet:
\begin{itemize}
\item Case 9: flux-type boundary condition (\ref{eq:inlet_flux_type}) with constant inlet concentration $c_{0}(t) = c_{0}$;
\item Case 10: flux-type boundary condition (\ref{eq:inlet_flux_type}) with step inlet concentration $c_{0}(t) = c_{0}H(t_{0}-t)$ and pulse duration $t_{0} = 3\,\textrm{days}$.
\end{itemize} 
In Figure \ref{fig:plots2}, we plot the relative concentration distributions ($c(x,t)/c_{0}$) obtained from our semi-analytical solution at several points in time. Only case 9 is considered by \citet{liu_1998} and \citet{guerrero_2013} who both report relative concentration distributions at $t = 2, 6, 10\,\textrm{days}$ only with their solutions at these times all in excellent agreement with those for case 9 in Figure \ref{fig:plots2}. For case 10 and the remaining times for case 9, agreement with the relative concentration distributions computed using the numerical method discussed in Appendix \ref{app:fvm} demonstrates that our semi-analytical solution produces the correct results.

To further highlight the capability of our approach, we consider three final test cases (cases 11-13 in Tables \ref{tab:transport_parameters} and \ref{tab:boundary_conditions}). Case 11 is case 9 with a continuously varying inlet concentration (see Table \ref{tab:boundary_conditions}). Case 12 considers the same five-layer medium as cases 9--11 but with a zero concentration gradient applied at both the inlet and outlet and an initial condition modelling injection of a contaminant into the medium at time $t = 0$ between $x = 14\,\mathrm{cm}$ and $x = 18\,\mathrm{cm}$ \citep{goltz_2017}. This test case is solved by creating an artificial layer extending from $x = 14\,\mathrm{cm}$ to $x = 18\,\mathrm{cm}$ having the same transport parameters as the sand layer in which it is located, leading to a seven-layer formulation as described in Table \ref{tab:transport_parameters}. Finally, case 13 demonstrates the full capability of our semi-analytical solution with non-zero rate constants of decay and production in each layer and multiple contaminant sources (see Tables \ref{tab:transport_parameters} and \ref{tab:boundary_conditions}). The agreement with the relative concentration distributions obtained using the numerical method evident in Figure \ref{fig:plots2} for cases 11--13 further confirm the correctness of our semi-analytical solution approach.

\section{Conclusion}
\label{sec:conclusions}
In this paper, we have developed a semi-analytical Laplace-transform based method to solve the one-dimensional linear advection-dispersion-reaction equation in a layered medium. The novelty of the approach is to introduce unknown functions at the interfaces between adjacent layers, which allows the multilayer problem to be isolated and solved separately on each layer before being numerically inverted back to the time domain. Our derived solution is quite general in that it can be applied to problems involving an arbitrary number of layers and arbitrary time-varying boundary conditions at the inlet and outlet. The derived solutions extend and generalise recent work on diffusion \citep{carr_2016,carr_2018a,rodrigo_2016} and reaction-diffusion \citep{zimmerman_2016} in layered media. 

The solutions presented in this paper, and our MATLAB code, are limited to constant initial conditions in each layer and interface conditions imposing continuity of concentration and dispersive flux between adjacent layers. However, extension to other interface conditions that do not impose continuity of concentration (see, e.g., \citet{carr_2018a}), is straightforward and achieved in our formulation by replacing Eq (\ref{eq:concentration_continuity_Laplace}) with the Laplace transform of the imposed condition. For example, solving the multilayer transport model (\ref{eq:pde})--(\ref{eq:general_bcR}) with concentration continuity (\ref{eq:concentration_continuity}) replaced by the \textit{partition interface condition} \citep{carr_2018a,carr_2018b}, $c_{i}(\ell_{i},t) = \alpha_{i}c_{i+1}(\ell_{i},t)$ where $\alpha_{i} > 0$ is a specified constant, simply requires replacement of Eq (\ref{eq:concentration_continuity_Laplace}) with $C_{i}(\ell_{i},s) = \alpha_{i}C_{i+1}(\ell_{i},s)$ and ultimately a small modification to the entries of the linear system (\ref{eq:linear_system}). Treatment of spatially-varying initial conditions, where $f_{i}$ is now $f_{i}(x)$ in the boundary value problems (\ref{eq:bvp1_de})--(\ref{eq:bvp1_bc2}), (\ref{eq:bvpi_de})--(\ref{eq:bvpi_bc2}) and (\ref{eq:bvpm_de})--(\ref{eq:bvpm_bc2}), is more challenging but possible and would lead to more complicated expressions for the Laplace domain concentration than those defined in Eqs (\ref{eq:C1})--(\ref{eq:Cm}) and Table \ref{tab:FAB}. {The method for inverting the Laplace transform can lead to unreliable results for advection-dominated transport so addressing this limitation would be valuable to pursue in the future.} Another possible avenue for future research is to extend the semi-analytical solutions to accommodate multispecies transport and/or rate-limited sorption, as described by \cite{chen_2019a,chen_2019b} for single-layer media and \cite{mieles_2012} and \cite{chen_2018} for two-layer media. We address the case of multispecies multilayer transport in a forthcoming paper \citep{carr_2020}.

\appendix
\section{Numerical solution}
\label{app:fvm}
In this appendix, we briefly outline the finite volume scheme used to obtain a numerical solution to the multilayer transport model (\ref{eq:pde})--(\ref{eq:general_bcR}), as mentioned in section \ref{sec:results}. 

The interval $[0,L]$ is discretised using a uniform grid consisting of $n$ nodes with the $k$th node located at $x = (k-1)h =: x_{k}$, where $k = 1,\hdots,n$ and $h = L/(n-1)$. The number of nodes $n$ is chosen to ensure that a node coincides with each interface ($x = \ell_{i}$, $i=1,\hdots,m-1$). Let $\overline{c}_{k}(t)$ denote the numerical approximation to $c(x,t)$ at $x = x_{k}$ and
\begin{align*}
J_{i,k} &= D_{i}\frac{\overline{c}_{k}-\overline{c}_{k-1}}{h} - v_{i}\frac{\overline{c}_{k-1}+\overline{c}_{k}}{2},\\
S_{i,k} &= -\mu_{i}\overline{c}_{k} + \gamma_{i}.
\end{align*}
The discrete system takes the form:
\begin{align}
\label{eq:ivp}
\mathbf{M}\frac{\textrm{d}\mathbf{c}}{\textrm{d}t} = \mathbf{F}(\mathbf{c}),\quad \mathbf{c}(0) = \mathbf{c}_{0},
\end{align}
where $\mathbf{M}\in\mathbb{R}^{n\times n}$, $\mathbf{c} = (\overline{c}_{1},\hdots,\overline{c}_{n})^{T}\in\mathbb{R}^{n}$ and $\mathbf{F} = (F_{1},\hdots,F_{n})^{T}\in\mathbb{R}^{n}$. The initial solution vector $\mathbf{c}_{0}\in\mathbb{R}^{n}$ gets its entries from the initial condition (\ref{eq:ic}) with first entry $f_{1}$, last entry $f_{m}$ and $k$th entry ($k = 2,\hdots,n-1$) equal to $f_{i}$ if $x_{k}\in(\ell_{i-1},\ell_{i})$ or $(f_{i}+f_{i+1})/2$ if $x_{k} = \ell_{i}$, where $i = 1,\hdots,m-1$ is the interface index. The form of $\mathbf{M}$ depends on the choice of boundary conditions at the inlet and outlet:
\begin{align*}
\mathbf{M} = \begin{cases} \mathbf{I}, & \text{if $b_{0}\neq 0$ and $b_{L}\neq 0$,}\\
\mathbf{I} - \mathbf{e}_{1}\mathbf{e}_{1}^{T}, & \text{if $b_{0} = 0$ and $b_{L}\neq 0$,}\\
\mathbf{I} - \mathbf{e}_{n}\mathbf{e}_{n}^{T}, & \text{if $b_{0} \neq 0$ and $b_{L} = 0$,}\\
\mathbf{I} - \mathbf{e}_{1}\mathbf{e}_{1}^{T} - \mathbf{e}_{n}\mathbf{e}_{n}^{T}, & \text{if $b_{0} = 0$ and $b_{L} = 0$,}\\ \end{cases}
\end{align*}
where $\mathbf{I}$ is the $n\times n$ identity matrix and $\mathbf{e}_{k}$ is the $k$th column of $\mathbf{I}$. The components of $\mathbf{F}$ are defined as follows:
\begin{align*}
F_{1} = a_{0}\overline{c}_{1} - g_{0}(t)
\end{align*}
if $b_{0} = 0$,
\begin{align*}
F_{1} = \frac{J_{1,2} + \frac{D_{1}}{b_{0}}g_{0}(t) + (v_{1} - D_{1}\frac{a_{0}}{b_{0}})\overline{c}_{1} + \frac{h}{2}S_{1,1}}{\frac{h}{2}R_{1}}
\end{align*}
if $b_{0}\neq 0$,
\begin{align*}
F_{k} = \frac{J_{i,k+1} - J_{i,k} + hS_{i,k}}{hR_{i}}
\end{align*}
if $x_{k}\in(\ell_{i-1},\ell_{i})$,
\begin{align*}
F_{k} = \frac{\theta_{i+1}J_{i+1,k+1} - \theta_{i}J_{i,k} + \frac{h}{2}\left(\theta_{i}S_{i,k} + \theta_{i+1}S_{i+1,k}\right)}{\frac{h}{2}\left(\theta_{i}R_{i} + \theta_{i+1}R_{i+1}\right)}
\end{align*}
if $x_{k} = \ell_{i}$ and $i = 2,\hdots,m-1$,
\begin{align*}
F_{n} = a_{L}\overline{c}_{n} - g_{L}(t)
\end{align*}
if $b_{L} = 0$, and
\begin{align*}
F_{n} = \frac{\frac{D_{m}}{b_{L}}g_{L}(t) - (v_{m} + D_{m}\frac{a_{L}}{b_{L}})\overline{c}_{n} - J_{m,n} + \frac{h}{2}S_{m,n}}{\frac{h}{2}R_{m}}
\end{align*}
if $b_{L}\neq 0$. The initial value problem (\ref{eq:ivp}) is solved using MATLAB's in-built \texttt{ode15s} solver with the default tolerances and options of $\texttt{Mass} = \mathbf{M}$ and $\texttt{MassSingular} = \texttt{true}$ if $b_{0}$ and/or $b_{L}$ are zero. The interval of integration (\texttt{tspan}) is chosen to return the solution at the appropriate times shown in Figures \ref{fig:plots1}--\ref{fig:plots2}.

\begin{acknowledgements}
This research was funded by the Australian Research Council (DE150101137). We thank the two anonymous reviewers for suggestions that improved the final manuscript.
\end{acknowledgements}

\bibliographystyle{spbasic}      
\bibliography{references}   

\end{document}